\documentclass[10pt]{amsart}
\usepackage{geometry} % see geometry.pdf on how to lay out the page. There's lots.
\usepackage{graphicx}
\usepackage{amssymb}
\usepackage{amsmath,amssymb,tikz-cd}
\usepackage{amsthm}
\usepackage{hyperref, enumitem, longtable, array}
\usepackage{mathrsfs}

\usepackage{setspace}%\doublespace 

\theoremstyle{plain}
\newtheorem{theorem}{Theorem}[section]
\newtheorem{proposition}[theorem]{Proposition}

\newtheorem{corollary}[theorem]{Corollary}

\theoremstyle{definition}

\theoremstyle{remark}
\newtheorem{remark}[theorem]{Remark}

\title{Universal Homotopy Theories and Associated Homological Algebras}
\author{Ahmad Rouintan}
\email{amanooman@gmail.com}
\date{} % delete this line to display the current date

\begin{document}
\setstretch{1.1}

\maketitle

\begin{abstract}
	Let $\mathscr{C}$ be a small category. For every commutative ring $R$ with unity, we associate an $R\mathrm{-linear}$ abelian category with the universal homotopy category of $\mathscr{C}$, where we can do the corresponding homological algebra. %Then, we apply this to the category of smooth schemes of finite type over a field $k$ of characteristic zero, to give an alternative description of Voevodsky's triangulated category of motives over $k$.
\end{abstract}

\tableofcontents

\emph{There was a problem in the final section of the previous version, so I deleted it and changed the title.}

\section*{Introduction}

%Let $Sm/k$ be the category of smooth schemes of finite type over a field $k$ of characteristic zero. To explain the $\mathbb{A}^1\mathrm{-homotopy}$ theory of Morel and Voevodsky, Daniel Dugger introduced the notion of universal homotopy theory for small categories in \cite{Dugger2001}. He proved that the $\mathbb{A}^1\mathrm{-homotopy}$ category can be recovered starting with the universal homotopy category of $Sm/k$ and performing two localizations. 

Let $\mathscr{M}$ be a model category. In Section 2.5 of \cite{Quillen1967}, Quillen states that under some conditions the category of commutative group objects $\mathscr{M}_{\mathbf{Ab}}\subseteq \mathscr{M}$ provides a framework to study (co)homology theories for $\mathscr{M}$. In other words, using $\mathscr{M}_{\mathbf{Ab}}$ we can \lq\lq linearize" or \lq\lq abelianize" the homotopical algebra of $\mathscr{M}$, resulting in a platform to do homological algebra. This short paper is a completion and generalization of this work of Quillen, although I did not know about it at the time I first wrote it. By the word completion, I mean that instead of working with a model category, we start with an arbitrary small category $\mathscr{C}$ and use its universal model category - a notion introduced by Daniel Dugger in \cite{Dugger2001}. This leads us to construct a category using commutative group objects $\mathscr{C}_{\mathbf{Ab}} \subseteq \mathscr{C}$ that satisfies all the conditions Quillen assumed and therefore, a platform to do homological algebra associated with the universal homotopy theory of $\mathscr{C}$. By the word generalization, I mean that we will develop our theory not only for commutative group objects but also for the more general category of $R\mathrm{-module}$ objects\footnote{This might not make sense in $\mathscr{C}$ itself, but we can make sense of it using the Yoneda embedding.} in $\mathscr{C}$, denoted by $\mathscr{C}_{R\mathbf{-Mod}}$. The idea, naively, is as follows:

\begin{quotation}
\emph{Let $R$ be a commutative ring with unity, and let $\mathscr{C}$ be a small category. Then, if we apply the same steps that construct the universal model category of $\mathscr{C}$ to the category of $R\mathit{-module}$ objects in $\mathscr{C}$, we get an $R\mathit{-linear}$ abelian category that is closely related to the universal homotopy category of $\mathscr{C}$.} 
\end{quotation}

%%%%%%%%%%%%%%%%%%%%
\section{Associated homological algebra with coefficients in $\mathbb{Z}$}\label{1}
%%%%%%%%%%%%%%%%%%%%

Let $\mathscr{C}$ be a small category with finite products and a terminal object. Then, we have the following proposition.

\begin{proposition}[\cite{MacLane1998}, Section 3.6, Proposition 1]\label{grpobj}
An object $g\in \mathscr{C}$ is a group object in $\mathscr{C}$ if and only if $\mathrm{Mor}_\mathscr{C}(-, g)$ is a group object in the functor category $\mathbf{Set}^{\mathscr{C}^\mathrm{op}}$.
\end{proposition}

Because the category $\mathbf{Set}^{\mathscr{C}^\mathrm{op}}$ has all products and a terminal object, we can omit these conditions on $\mathscr{C}$ and call an object $g\in \mathscr{C}$ a group object, whenever $\mathrm{Mor}_\mathscr{C}(-, g)$ is a group object in $\mathbf{Set}^{\mathscr{C}^\mathrm{op}}$. Similarly, we call an object $g\in \mathscr{C}$ a commutative group object, whenever $\mathrm{Mor}_\mathscr{C}(-, g)$ is a commutative group object in $\mathbf{Set}^{\mathscr{C}^\mathrm{op}}$. A morphism $f\in \mathrm{Mor}_{\mathscr{C}}(g, h)$  is a morphism of commutative group objects the Yoneda embedding takes it to a morphism of commutative group objects from $\mathrm{Mor}_\mathscr{C}(-, g)$ to $\mathrm{Mor}_\mathscr{C}(-, h)$. We denoted the subcategory of commutative group object in $\mathscr{C}$ by $\mathscr{C}_{\mathbf{Ab}}$. 

\begin{proposition}\label{abobj}
Commutative group objects in $\mathbf{Set}^{\mathscr{C}^\mathrm{op}}$ are exactly presheaves of abelian groups on $\mathscr{C}$. Moreover, the Yoneda embedding identifies the subcategory of commutative group objects $\mathscr{C}_{\mathbf{Ab}} \subseteq \mathscr{C}$ with representable presheaves of abelian groups on $\mathscr{C}$.
\end{proposition}

\begin{proof}
Let $\mathcal{G}$ be a commutative group object in $\mathbf{Set}^{\mathscr{C}^\mathrm{op}}$. Then, if we consider the diagrams that give it the structure of a commutative group object and apply those diagrams to each object $c\in\mathscr{C}$, we see that $\mathcal{G}$ is a presheaf of abelian groups on $\mathscr{C}$. On the other hand, if $\mathcal{G}$ is a presheaf of abelian groups on $\mathscr{C}$, then, $\mathcal{G}(c)$ for each object $c\in \mathscr{C}$ is an abelian group so we have a multiplication $\mu_c: \mathcal{G}(c) \times \mathcal{G}(c)\to \mathcal{G}(c)$, an identity $\iota_c: 0\to \mathcal{G}(c)$, and an inverse map $\eta_c: \mathcal{G}(c)\to \mathcal{G}(c)$ that satisfy the conditions for $\mathcal{G}(c)$ to be an abelian group. Combining all of these, we get a multiplication map $\mu: \mathcal{G}\times \mathcal{G}\to \mathcal{G}$, an identity map $\iota: 0\to \mathcal{G}$ where $0$ is the constant presheaf, and an inverse map $\eta: \mathcal{G}\to \mathcal{G}$ satisfying the conditions for $\mathcal{G}$ to be a commutative group object in $\mathbf{Set}^{\mathscr{C}^\mathrm{op}}$.\\
For the second part, notice that by definition the embedding of $\mathscr{C}_{\mathbf{Ab}}$ into the subcategory of representable presheaves of abelian groups is fully faithful. It also is surjective on objects and that completes the proof.
\end{proof}

We have a free-forgetful adjunction between $\mathbf{Set}^{\mathscr{C}^\mathrm{op}}$ and $\mathbf{Ab}^{\mathscr{C}^\mathrm{op}}$ by applying the usual free-forgetful adjunction between sets and abelian groups, object-wise. We denote both of these adjunctions by $f_\mathbb{Z} \dashv U$.  All of the above discussions can be summarised in the diagram

\begin{equation}\label{cocom}
	\begin{tikzcd}
		\mathscr{C}
		\arrow[hookrightarrow, r, "Y"] &
		\mathbf{Set}^{\mathscr{C}^\mathrm{op}}
		\arrow[d, shift right=0.5ex, swap, "{f_\mathbb{Z}}"] \\
		\mathscr{C}_{\mathbf{Ab}}
		\arrow[hookrightarrow, u]
		\arrow[hookrightarrow, r, "Y"] &
		\mathbf{Ab}^{\mathscr{C}^\mathrm{op}}
		\arrow[u, shift right=0.5ex, swap, "U"]
	\end{tikzcd}
\end{equation} 
where $Y$ is the Yoneda embedding. This diagram is commutative if we only consider the forgetful functor $U$ from the free-forgetful adjunction.

\begin{remark}
The constant presheaf $\mathbb{Z}$ is a commutative ring object in $\mathbf{Set}^{\mathscr{C}^\mathrm{op}}$, so we can form its subcategory of $\mathbb{Z}\mathrm{-modules}$ which coincides with $\mathbf{Ab}^{\mathscr{C}^\mathrm{op}}$. The free-forgetful adjunction in Diagram \ref{cocom} just forgets the $\mathbb{Z}\mathrm{-module}$ structure in one direction and constructs the free $\mathbb{Z}\text{-module}$ in the other. This observation is crucial in the later parts of this paper.
\end{remark}

Here, let's take a moment to explain why we are interested in $\mathbf{Set}^{\mathscr{C}^\mathrm{op}}$ and $\mathbf{Ab}^{\mathscr{C}^\mathrm{op}}$. Consider the category of sets, the category of abelian groups, and the free-forgetful adjunction between them. If we take simplicial objects built out of sets, we get the category of simplicial sets, which is the home for algebraic topology. On the other hand, if we take simplicial objects built out of abelian groups, we get the category of simplicial abelian groups, which is the home to the classical homological algebra because it is equivalent to the category of connective chain complexes of abelian groups. Also, the free-forgetful adjunction has a simplicial version that links algebraic topology to homological algebra. By considering the two functor categories in Diagram \ref{cocom}, which inherit properties of sets and abelian groups respectively, we are taking the first baby step toward constructing a homotopy theory for $\mathscr{C}$ and associating a homological algebra with it. 

Fortunately, a bigger step toward the homotopy theory of $\mathscr{C}$ has already been taken by Daniel Dugger in \cite{Dugger2001}. First, let's consider the category of simplicial objects built out of objects of $\mathbf{Set}^{\mathscr{C}^\mathrm{op}}$ i.e. $\Delta^{\text{op}}\mathbf{Set}^{\mathscr{C}^\mathrm{op}}$. Using the tensor-Hom adjunction and the symmetry of the closed monoidal structure on the category of categories, this is the same as the category of presheaves of simplicial sets on $\mathscr{C}$ i.e. $s\mathbf{Set}^{\mathscr{C}^\mathrm{op}}$. As shown in \cite{Dugger2001}, this category has an object-wise model structure that turns it into the universal model category of $\mathscr{C}$. We introduce this model structure through the next theorem and corollary.

\begin{theorem}[Quillen, \cite{Quillen1967}, II, \S 2, Theorem 2.4]\label{Qui}
Let $\mathscr{D}$ be a category closed under finite limits and having sufficient many projectives. Let $s\mathscr{D}$ be the simplicial category of simplicial objects over $\mathscr{D}$. Define a map $\phi$ in $s\mathscr{D}$ to be a fibration (resp. weak equivalence) if $\mathrm{Hom}(P,\phi)$ is a fibration (resp. weak equivalence) in $s\mathbf{Set}$ for each projective object $P$ of $\mathscr{D}$, and a cofibration if $\phi$ has the Left lifting property with respect to the class of trivial fibrations. Then $s\mathscr{D}$ is a closed simplicial model category if $\mathscr{D}$ satisfies one of the following extra conditions:
\begin{enumerate}
\item
Every object of $s\mathscr{D}$ is fibrant,
\item
$\mathscr{D}$ is closed under inductive limits and has a set of small projective generators.
\end{enumerate}

\end{theorem}

\begin{corollary}[\cite{BousfieldKan1972}, Page 314]\label{BKmodel}
Define a morphism $\phi: \mathcal{O}\to \mathcal{O}' \in s\mathbf{Set}^{\mathscr{C}^\mathrm{op}}$ to be an object-wise weak equivalence (resp. object-wise fibration) if for every object $c\in \mathscr{C}$ the induced map $\phi_c: \mathcal{O}(c)\to \mathcal{O}'(c)$ is a weak equivalence (resp. fibration) of simplicial sets. These classes of object-wise weak equivalences and object-wise fibrations define a closed simplicial model structure on $s\mathbf{Set}^{\mathscr{C}^\mathrm{op}}$. 
\end{corollary}

The category $\mathbf{Set}^{\mathscr{C}^\mathrm{op}}$ with the model structure of Corollary \ref{BKmodel}, which is called the Bousfield-Kan model structure, turns out to be the universal model category of $\mathscr{C}$, in the sense of the following proposition. 

\begin{proposition}[\cite{Dugger2001}, Proposition 2.3]\label{univhom}
Any functor $\mathscr{F}: \mathscr{C} \to \mathscr{M}$ into a model category $\mathscr{M}$ may be factored through $s\mathbf{Set}^{\mathscr{C}^\mathrm{op}}$, in the sense that there is a Quillen pair 
\begin{equation*}
\begin{tikzcd}
\mathrm{Re}_{\mathscr{F}}: s\mathbf{Set}^{\mathscr{C}^\mathrm{op}}
\arrow[r, shift left=0.5ex] &
\mathscr{M}: \mathrm{Sing}_{\mathscr{F}}
\arrow[l, shift left=0.5ex]
\end{tikzcd}
\end{equation*} 
where $s\mathbf{Set}^{\mathscr{C}^\mathrm{op}}$ is considered with the Bousfield-Kan model structure and the diagram
\begin{equation*}
\begin{tikzcd}
\mathscr{C} 
\arrow[hookrightarrow, r] 
\arrow[rd, swap, "{\mathscr{F}}"] & 
s\mathbf{Set}^{\mathscr{C}^\mathrm{op}}
\arrow[d, "{\mathrm{Re}_{\mathscr{F}}}"] \\
& 
\mathscr{M}
\end{tikzcd}
\end{equation*}
commutes up to a natural weak equivalence. Moreover, the category of such factorizations is contractible. 
\end{proposition}

We denote the homotopy category of $s\mathbf{Set}^{\mathscr{C}^\mathrm{op}}$ with respect to the Bousfield-Kan model structure by $\mathscr{H}(\mathscr{C})$. This is not an ideal notation but makes sense because it is the universal homotopy category of $\mathscr{C}$.

Having Theorem \ref{univhom}, we want to use the free-forgetful adjunction to construct a corresponding category where we can do homological algebra. First, we need to consider the category of simplicial objects built out of objects of $\mathbf{Ab}^{\mathscr{C}^\mathrm{op}}$ i.e. $\Delta^{\text{op}}\mathbf{Ab}^{\mathscr{C}^\mathrm{op}}$. This is again the same as the category of presheaves on $\mathscr{C}$ with values in $s\mathbf{Ab}$ i.e. $s\mathbf{Ab}^{\mathscr{C}^\mathrm{op}}$. Here, the Dold-Kan correspondence gives us an equivalence between $\Delta^{\text{op}}\mathbf{Ab}^{\mathscr{C}^\mathrm{op}}$ and the category of connective chain complexes of $\mathbf{Ab}^{\mathscr{C}^\mathrm{op}}$, denoted by $C_+(\mathbf{Ab}^{\mathscr{C}^\mathrm{op}})$.

We expect that because the Bousfield-Kan model structure on $s\mathbf{Set}^{\mathscr{C}^\mathrm{op}}$ was defined object-wise, it induces a model structure on its subcategory $s\mathbf{Ab}^{\mathscr{C}^\mathrm{op}}$ through the free-forgetful adjunction, just like the case of $s\mathbf{Set}$ and $s\mathbf{Ab}$. But, we should find a set of projective generators for $\mathbf{Ab}^{\mathscr{C}^\mathrm{op}}$ in order to use Theorem \ref{Qui} to prove this claim. In the case of $\mathbf{Set}^{\mathscr{C}^\mathrm{op}}$, the set of small projective generators $\{\mathcal{P}_a\}_{a\in \mathscr{C}}$ is characterized by the natural isomorphisms 
$$\mathrm{Mor}_{\mathbf{Set}^{\mathscr{C}^\mathrm{op}}}(\mathcal{P}_a, \mathcal{B}) = \mathcal{B}(a)$$
for every presheaf $\mathcal{B}$. For $\mathbf{Ab}^{\mathscr{C}^\mathrm{op}}$, the set of small projective generators is given by $\{f_{\mathbb{Z}}(\mathcal{P}_a)\}_{a\in \mathscr{C}}$. Notice that For every presheaf of abelian groups $\mathcal{G}$ on $\mathscr{C}$ we have 
$$\mathcal{G}(a) = \mathrm{Hom}_{\mathbf{Ab}}(\mathbb{Z}, \mathcal{G}(a)) \cong  \mathrm{Hom}_{\mathbf{Ab}^{\mathscr{C}^\mathrm{op}}}(f_\mathbb{Z}(\mathcal{P}_a), \mathcal{G}).$$ %Check the proof!

Before introducing the model structure on $s\mathbf{Ab}^{\mathscr{C}^\mathrm{op}}$, recall that $s\mathbf{Ab}$ admits a model structure in which weak equivalences (resp. fibrations) are weak equivalences (resp. fibrations) of underlying simplicial sets. With this model structure, the free-forgetful adjunction becomes a Quillen adjunction. Using this, we can get an object-wise model structure on $s\mathbf{Ab}^{\mathscr{C}^\mathrm{op}}$ such that the corresponding free-forgetful adjunction becomes a Quillen adjunction.

\begin{corollary}\label{aAbpremodel}
Define a morphism $\phi: \mathcal{G}\to \mathcal{G}' \in s\mathbf{Ab}^{\mathscr{C}^\mathrm{op}}$ to be an object-wise weak equivalence (resp. object-wise fibration) if for every object $c\in \mathscr{C}$ the induced map $\phi_c: \mathcal{G}(c)\to \mathcal{G}'(c)$ is a weak equivalence (resp. fibration) of underlying simplicial sets. These classes of object-wise weak equivalences and object-wise fibrations define a closed simplicial model structure on $s\mathbf{Ab}^{\mathscr{C}^\mathrm{op}}$. Moreover, this model structure turns the free-forgetful adjunction between $s\mathbf{Set}^{\mathscr{C}^\mathrm{op}}$ and $s\mathbf{Ab}^{\mathscr{C}^\mathrm{op}}$ into a Quillen adjunction.
\end{corollary}

\begin{proof}
The first part of the proof is similar to the proof of Corollary \ref{BKmodel} with $\mathbf{Ab}^{\mathscr{C}^\mathrm{op}}$ instead of $\mathbf{Set}^{\mathscr{C}^\mathrm{op}}$. So, we only need to prove that $f_\mathbb{Z} \dashv U$ is a Quillen adjunction. This follows from the fact that the usual free-forgetful adjunction is a Quillen adjunction between $s\mathbf{Set}$ and $s\mathbf{Ab}$ and the fact that the model structures on $s\mathbf{Set}^{\mathscr{C}^\mathrm{op}}$ and $s\mathbf{Ab}^{\mathscr{C}^\mathrm{op}}$ are defined object-wise.
\end{proof}

Having the model structure on $\Delta^{\text{op}}\mathbf{Ab}^{\mathscr{C}^\mathrm{op}}$ we can use the Dold-Kan equivalence to impose a model structure on $C_+(\mathbf{Ab}^{\mathscr{C}^\mathrm{op}})$. This would turn weak equivalences into quasi-isomorphisms just like the usual Dold-Kan correspondence because we have defined all of the model structures object-wise. We denote the respective homotopy category of $C_+(\mathbf{Ab}^{\mathscr{C}^\mathrm{op}})$ by $\mathscr{H}_+(\mathscr{C}_{\mathbf{Ab}})$, which again is not an ideal notation! Lastly, because the free-forgetful adjunction is a Quillen adjunction based on Corollary \ref{aAbpremodel}, we get an adjoint pair between $\mathscr{H}(\mathscr{C})$ and $\mathscr{H}_+(\mathscr{C}_{\mathbf{Ab}})$ by taking the derived functors. We denote this adjunction by $f_\mathbb{Z} \dashv U$, too. 

We summarise all of these in the diagram

\begin{equation}\label{unsuniv}
	\begin{tikzcd}
		\mathscr{C}
		\arrow[hookrightarrow, r, "Y"] &
		\mathbf{Set}^{\mathscr{C}^\mathrm{op}}
		\arrow[d, shift right=0.5ex, swap, "{f_\mathbb{Z}}"] 
		\arrow[r] &
		\Delta^{\text{op}}\mathbf{Set}^{\mathscr{C}^\mathrm{op}}
		\arrow[r]
		\arrow[d, shift right=0.5ex, swap, "{f_\mathbb{Z}}"] &
		\mathscr{H}(\mathscr{C})
		\arrow[d, shift right=0.5ex, swap, "{f_\mathbb{Z}}"] \\
		\mathscr{C}_{\mathbf{Ab}}
		\arrow[hookrightarrow, u]
		\arrow[hookrightarrow, r, "Y"] &
		\mathbf{Ab}^{\mathscr{C}^\mathrm{op}}
		\arrow[r]
		\arrow[u, shift right=0.5ex, swap, "U"] &
		C_+(\mathbf{Ab}^{\mathscr{C}^\mathrm{op}}) 
		\arrow[r]
		\arrow[u, shift right=0.5ex, swap, "U"] &
		\mathscr{H}_+(\mathscr{C}_{\text{Ab}}).
		\arrow[u, shift right=0.5ex, swap, "U"]
	\end{tikzcd}
\end{equation} 
which is commutative if we only consider $U$ from the adjunction $f_\mathbb{Z} \dashv U$ in every column.

%\begin{remark}
%After writing the first draft of this paper, I realized that Quillen had discussed similar things in Section 2.5 of \cite{Quillen1967} where he assumes a few conditions on $\mathscr{C}_{\mathrm{Ab}}$ to study homology and cohomology theories. He even proves that for an abelian category $\mathscr{A}$ with enough projectives, the category of simplicial objects built out of it, denoted by $s\mathscr{A}$ satisfies those conditions. However, he could not take the final step because he was not aware of the notion of universal homotopy theories - a step that has been taken in this paper without prior knowledge of Quillen's work. For interested readers, it would be a great experience to go through Quillen's work while reading this paper.
%\end{remark}

Now that we have the homotopy categories, or better said unstable homotopy categories, it is straightforward to construct the stable homotopy categories. Let's focus on $\mathscr{H}_+(\mathscr{C}_{\mathbf{Ab}})$ first. It is not a triangulated category with respect to the shift functor. But, an obvious way to address this problem is to extend our category and consider all the chain complexes of $\mathbf{Ab}^{\mathscr{C}^\mathrm{op}}$ instead of only connective ones from the beginning. Then, by inverting quasi-isomorphisms, we get the unbounded derived category of chain complexes of $\mathbf{Ab}^{\mathscr{C}^\mathrm{op}}$, which we denote ignorantly by $D(\mathscr{C}_{\mathbf{Ab}})$. We will not get into more details here because the construction of $D(\mathscr{C}_{\mathbf{Ab}})$ is a classic that can be found in any reference on homological algebra - see Chapter 10 of \cite{Weibel1994}. 

On the other hand, to construct the stable homotopy category of $s\mathbf{Set}^{\mathscr{C}^\mathrm{op}}$ we need to stabilize with respect to the constant presheaf $S^1$. This has been done by John Frederic Jardine in \cite{Jardines1987} so we will not dwell on it further. One thing to consider here is the fact that in this stage we need to work with the pointed presheaves of simplicial sets on $\mathscr{C}$. However, this is not a problem as the free-forgetful functor we were working with factors through the free-forgetful functor $f_+ \dashv U$ between the unpointed and pointed settings. Also, note that $s\mathbf{Set}^{\mathscr{C}^\mathrm{op}}_\bullet$ is the universal pointed model category of $\mathscr{C}$ - see page 160 in \cite{Dugger2001}.

All of these can be summarised in the diagram

\begin{equation}\label{premap}
	\begin{tikzcd}
		\mathscr{C}
		\arrow[hookrightarrow, r, "Y"] &
		\mathbf{Set}^{\mathscr{C}^\mathrm{op}}
		\arrow[d, shift right=0.5ex, swap, "{f_+}"] 
		\arrow[r] &
		\Delta^{\text{op}}\mathbf{Set}^{\mathscr{C}^\mathrm{op}}
		\arrow[r]
		\arrow[d, shift right=0.5ex, swap, "{f_+}"] &
		\mathscr{H}(\mathscr{C})  
		\arrow[d, shift right=0.5ex, swap, "{f_+}"] &
		& \\
		&
		\mathbf{Set}^{\mathscr{C}^\mathrm{op}}_\bullet
		\arrow[d, shift right=0.5ex, swap, "{f_\mathbb{Z}}"]
		\arrow[u, shift right=0.5ex, swap, "U"] 
		\arrow[r] &
		\Delta^{\text{op}}\mathbf{Set}^{\mathscr{C}^\mathrm{op}}_\bullet
		\arrow[r]
		\arrow[d, shift right=0.5ex, swap, "{f_\mathbb{Z}}"]
		\arrow[u, shift right=0.5ex, swap, "U"] &
		\mathscr{H}_\bullet(\mathscr{C}) 
		\arrow[r, shift left=0.5ex, "{\Sigma_{S^1}^\infty}"] 
		\arrow[d, shift right=0.5ex, swap, "{f_\mathbb{Z}}"]
		\arrow[u, shift right=0.5ex, swap, "U"] &
		\mathscr{SH}(\mathscr{C})
		\arrow[l, shift left=0.5ex, "{\Omega_{S^1}^\infty}"]
		\arrow[d, shift right=0.5ex, swap, "{f_\mathbb{Z}}"] \\
		\mathscr{C}_{\mathbf{Ab}}
		\arrow[hookrightarrow, uu]
		\arrow[hookrightarrow, r, "Y"] &
		\mathbf{Ab}^{\mathscr{C}^\mathrm{op}}
		\arrow[r]
		\arrow[u, shift right=0.5ex, swap, "U"] &
		C_+(\mathbf{Ab}^{\mathscr{C}^\mathrm{op}}) 
		\arrow[r]
		\arrow[u, shift right=0.5ex, swap, "U"] &
		\mathscr{H}_+(\mathscr{C}_{\text{Ab}})
		\arrow[u, shift right=0.5ex, swap, "U"] 
		\arrow[r, shift left=0.5ex, "{\Sigma^\infty}"] &
		D(\mathscr{C}_{\text{Ab}})
		\arrow[l, shift left=0.5ex, "{\Omega^\infty}"]
		\arrow[u, shift right=0.5ex, swap, "U"]
	\end{tikzcd}
\end{equation} 
where $\mathscr{SH}(\mathscr{C})$ is the stable homotopy category of $s\mathbf{Set}^{\mathscr{C}^\mathrm{op}}_\bullet$. Again, this diagram becomes commutative if we only consider forgetful functors in every free-forgetful adjunction in the diagram.

%%%%%%%%%%%%%%%%%%%%
\section{Associated homological algebra with coefficients in an arbitrary ring $R$}\label{2}
%%%%%%%%%%%%%%%%%%%%

Let $\mathscr{C}$ be a small category. We will go back to the beginning and take similar steps, this time for an arbitrary commutative ring $R$ instead of $\mathbb{Z}$. First, notice that the constant presheaf $R$ in $\mathbf{Set}^{\mathscr{C}^\mathrm{op}}$ is a commutative ring object where $\mathscr{C}$ is a small category. Therefore, we can construct the subcategory of $R\mathrm{-modules}$ in $\mathbf{Set}^{\mathscr{C}^\mathrm{op}}$. Similar to the definition of commutative group objects in $\mathscr{C}$, we define an object $m\in \mathscr{C}$ to be an $R\mathrm{-module}$ object whenever $\mathrm{Mor}_\mathscr{C}(-, m)$ is an $R\mathrm{-module}$ object in $\mathbf{Set}^{\mathscr{C}^\mathrm{op}}$. We can state the analogue of Proposition \ref{abobj} for $R\mathrm{-module}$ objects. 

\begin{proposition}\label{Rmodobj}
The $R\mathit{-module}$ objects in $\mathbf{Set}^{\mathscr{C}^\mathrm{op}}$ are exactly presheaves of $R\mathit{-modules}$ on $\mathscr{C}$. Moreover, the Yoneda embedding identifies the subcategory of $R\mathit{-module}$ objects $\mathscr{C}_{R\mathbf{-Mod}} \subseteq \mathscr{C}$ with representable presheaves of $R\mathit{-modules}$ on $\mathscr{C}$.
\end{proposition}

\begin{proof}
An argument similar to Proposition \ref{abobj} gives the result.
\end{proof}

Obviously, the free-forgetful adjunction $f_R \dashv U$ between sets and $R\mathrm{-modules}$ has a presheaf version which we again denote by $f_R \dashv U$. We can summarize the above discussions in the diagram
\begin{equation}\label{cocomRmod}
	\begin{tikzcd}
		\mathscr{C}
		\arrow[hookrightarrow, r, "Y"] &
		\mathbf{Set}^{\mathscr{C}^\mathrm{op}}
		\arrow[d, shift right=0.5ex, swap, "{f_R}"] \\
		\mathscr{C}_{R\mathbf{-Mod}}
		\arrow[hookrightarrow, u]
		\arrow[hookrightarrow, r, "Y"] &
		R\mathbf{-Mod}^{\mathscr{C}^\mathrm{op}}
		\arrow[u, shift right=0.5ex, swap, "U"]
	\end{tikzcd}
\end{equation} 
where $Y$ is the Yoneda embedding. This diagram is commutative if we only consider the forgetful functor U from the free-forgetful adjunction.

Next, we need to consider the category of simplicial objects built out of objects in $R\mathbf{-Mod}^{\mathscr{C}^\mathrm{op}}$ i.e.  $sR\mathbf{-Mod}^{\mathscr{C}^\mathrm{op}} \cong \Delta^{\mathrm{op}}R\mathbf{-Mod}^{\mathscr{C}^\mathrm{op}}$, and define a model structure on it. To do so, let's first recall the model structure for the category of simplicial $R\mathrm{-modules}$. The category $sR\mathbf{-Mod}$ admits a model structure similar to the model structure on $s\mathbf{Ab}$ in which weak equivalences (resp. fibrations) where weak equivalences (resp. fibrations) of underlying simplicial sets. On the other hand, recall that if we start with a simplicial $R\mathrm{-module}$ and construct its normalized chain complex, we get a connective chain complex of $R\mathrm{-modules}$. This defines a functor which is an equivalence of categories between the category of simplicial $R\mathrm{-module}$ $sR\mathbf{-Mod}$ and the category of connective chain complexes of $R\mathrm{-modules}$ $C_+(R\mathbf{-Mod})$. It is basically the $R\mathrm{-module}$ version of the Dold-Kan correspondence. So, we transfer the model structure of $sR\mathbf{-Mod}$ to $C_+(R\mathbf{-Mod})$ using this equivalence of categories. Now, we have all the tools needed to define an object-wise model structure on $sR\mathbf{-Mod}^{\mathscr{C}^\mathrm{op}}$.

\begin{corollary}\label{aRmodpremodel}
Define a morphism $\phi: \mathcal{G}\to \mathcal{G}' \in sR\mathbf{-Mod}^{\mathscr{C}^\mathrm{op}}$ to be an object-wise weak equivalence (resp. object-wise fibration) if for every object $c\in \mathscr{C}$ the induced map $\phi_c: \mathcal{G}(c)\to \mathcal{G}'(c)$ is a weak equivalence (resp. fibration) of underlying simplicial sets. These classes of object-wise weak equivalences and object-wise fibrations define a closed simplicial model structure on $sR\mathbf{-Mod}^{\mathscr{C}^\mathrm{op}}$. Moreover, this model structure turns the free-forgetful adjunction into a Quillen adjunction between $s\mathbf{Set}^{\mathscr{C}^\mathrm{op}}$ and $sR\mathbf{-Mod}^{\mathscr{C}^\mathrm{op}}$.
\end{corollary}

\begin{proof}
The first part of the proof is similar to the proof of Corollary \ref{BKmodel} with $R\mathbf{-Mod}^{\mathscr{C}^\mathrm{op}}$ instead of $\mathbf{Set}^{\mathscr{C}^\mathrm{op}}$. So, we only need to prove that $f_R \dashv U$ is a Quillen adjunction. This follows from the fact that the usual free-forgetful adjunction is a Quillen adjunction between $s\mathbf{Set}$ and $sR\mathbf{-Mod}$ and the fact that the model structures on $s\mathbf{Set}^{\mathscr{C}^\mathrm{op}}$ and $sR\mathbf{-Mod}^{\mathscr{C}^\mathrm{op}}$ are defined object-wise.
\end{proof}

This model structure can be transferred to $C_+(R\mathbf{-Mod}^{\mathscr{C}^\mathrm{op}})$ by the Dold-Kan equivalence, and its weak equivalences would be quasi-isomorphisms. We denote the respective homotopy category by $\mathscr{H}_+(\mathscr{C}_{R\mathbf{-Mod}})$. Therefore, we can summarize the above discussions in the diagram

\begin{equation}\label{unsunivRmod}
	\begin{tikzcd}
		\mathscr{C}
		\arrow[hookrightarrow, r, "Y"] &
		\mathbf{Set}^{\mathscr{C}^\mathrm{op}}
		\arrow[d, shift right=0.5ex, swap, "{f_R}"] 
		\arrow[r] &
		\Delta^{\text{op}}\mathbf{Set}^{\mathscr{C}^\mathrm{op}}
		\arrow[r]
		\arrow[d, shift right=0.5ex, swap, "{f_R}"] &
		\mathscr{H}(\mathscr{C})
		\arrow[d, shift right=0.5ex, swap, "{f_R}"] \\
		\mathscr{C}_{R\mathbf{-Mod}}
		\arrow[hookrightarrow, u]
		\arrow[hookrightarrow, r, "Y"] &
		R\mathbf{-Mod}^{\mathscr{C}^\mathrm{op}}
		\arrow[r]
		\arrow[u, shift right=0.5ex, swap, "U"] &
		C_+(R\mathbf{-Mod}^{\mathscr{C}^\mathrm{op}}) 
		\arrow[r]
		\arrow[u, shift right=0.5ex, swap, "U"] &
		\mathscr{H}_+(\mathscr{C}_{R\mathbf{-Mod}})
		\arrow[u, shift right=0.5ex, swap, "U"]
	\end{tikzcd}
\end{equation} 
which is commutative if we only consider $U$ from the adjunction $f_\mathbb{Z} \dashv U$ in every column. This diagram is analogous to Diagram \ref{unsuniv}. 

Lastly, we can address the fact that $\mathscr{H}_+(\mathscr{C}_{R\mathbf{-Mod}})$ is not triangulated by extending to the unbounded derived category $D(\mathscr{C}_{R\mathbf{-Mod}})$. Storing all of this information in a single diagram brings us to the analogue of Diagram \ref{premap} which looks like the following:

\begin{equation}\label{premapRmod}
	\begin{tikzcd}
		\mathscr{C}
		\arrow[hookrightarrow, r, "Y"] &
		\mathbf{Set}^{\mathscr{C}^\mathrm{op}}
		\arrow[d, shift right=0.5ex, swap, "{f_+}"] 
		\arrow[r] &
		\Delta^{\text{op}}\mathbf{Set}^{\mathscr{C}^\mathrm{op}}
		\arrow[r]
		\arrow[d, shift right=0.5ex, swap, "{f_+}"] &
		\mathscr{H}(\mathscr{C})  
		\arrow[d, shift right=0.5ex, swap, "{f_+}"] &
		& \\
		&
		\mathbf{Set}^{\mathscr{C}^\mathrm{op}}_\bullet
		\arrow[d, shift right=0.5ex, swap, "{f_R}"]
		\arrow[u, shift right=0.5ex, swap, "U"] 
		\arrow[r] &
		\Delta^{\text{op}}\mathbf{Set}^{\mathscr{C}^\mathrm{op}}_\bullet
		\arrow[r]
		\arrow[d, shift right=0.5ex, swap, "{f_R}"]
		\arrow[u, shift right=0.5ex, swap, "U"] &
		\mathscr{H}_\bullet(\mathscr{C}) 
		\arrow[r, shift left=0.5ex, "{\Sigma_{S^1}^\infty}"] 
		\arrow[d, shift right=0.5ex, swap, "{f_R}"]
		\arrow[u, shift right=0.5ex, swap, "U"] &
		\mathscr{SH}(\mathscr{C})
		\arrow[l, shift left=0.5ex, "{\Omega_{S^1}^\infty}"]
		\arrow[d, shift right=0.5ex, swap, "{f_R}"] \\
		\mathscr{C}_{R\mathbf{-Mod}}
		\arrow[hookrightarrow, uu]
		\arrow[hookrightarrow, r, "Y"] &
		R\mathbf{-Mod}^{\mathscr{C}^\mathrm{op}}
		\arrow[r]
		\arrow[u, shift right=0.5ex, swap, "U"] &
		C_+(R\mathbf{-Mod}^{\mathscr{C}^\mathrm{op}}) 
		\arrow[r]
		\arrow[u, shift right=0.5ex, swap, "U"] &
		\mathscr{H}_+(\mathscr{C}_{R\mathbf{-Mod}})
		\arrow[u, shift right=0.5ex, swap, "U"] 
		\arrow[r, shift left=0.5ex, "{\Sigma^\infty}"] &
		D(\mathscr{C}_{R\mathbf{-Mod}}).
		\arrow[l, shift left=0.5ex, "{\Omega^\infty}"]
		\arrow[u, shift right=0.5ex, swap, "U"]
	\end{tikzcd}
\end{equation} 
Again, this diagram is commutative if we only consider forgetful functors in every free-forgetful adjunction in the diagram.

\begin{remark}
In this section, we worked with modules over commutative ring objects that were constant as a presheaf. But, the constructions in this section can be done for arbitrary commutative ring objects in $\mathbf{Set}^{\mathscr{C}^\mathrm{op}}$. Moreover, the same can be done for any $R\mathrm{-linear}$ abelian model category $\mathscr{M}$ instead of $C_+(R\mathbf{-Mod}^{\mathscr{C}^\mathrm{op}})$ and a Quillen adjunction 
\begin{equation*}
\begin{tikzcd}
\mathscr{L}: s\mathbf{Set}^{\mathscr{C}^\mathrm{op}}
\arrow[r, shift left=0.5ex] &
\mathscr{M}: \mathscr{R}.
\arrow[l, shift left=0.5ex]
\end{tikzcd}
\end{equation*} 
\end{remark}

Lastly, the categories $\mathbf{Ab}^{\mathscr{C}^\mathrm{op}}$ and $R\mathbf{-Mod}^{\mathscr{C}^\mathrm{op}}$ are abelian categories with enough projectives, therefore, $s\mathbf{Ab}^{\mathscr{C}^\mathrm{op}}$ and $sR\mathbf{-Mod}^{\mathscr{C}^\mathrm{op}}$ satisfy the conditions stated in Section 2.5 of \cite{Quillen1967} by Quillen. This completes our task.

%\input{Section2.tex}

%\input{Section3.tex}

%-----------------------------------------------------------------------------------

\end{document}